\documentclass[11pt,reqno]{amsart} % please use amsart at 11pt

\usepackage{amssymb,latexsym}
\usepackage{cite} % to get Refs. [1,2,3] typeset as [1--3] automatically

\usepackage{fontenc}
\usepackage{times}

\usepackage{amsfonts,amssymb,epsfig}
\usepackage{amscd}
\usepackage{color}
 \usepackage{hyperref}

\usepackage{textcomp}
\usepackage{mathrsfs}

\usepackage[height=200mm,width=140mm]{geometry} % this is the journal
\theoremstyle{plain}

\theoremstyle{remark}

\theoremstyle{example}

\numberwithin{equation}{section} % to get equations numbered
\theoremstyle{definition}
\newtheorem{definition}{Definition}[section]

\newtheorem{example}[definition]{Example}
\newtheorem{remark}[definition]{Remark}

\newtheorem{theorem}[definition]{Theorem}
\newtheorem{proposition}[definition]{Proposition}

\newtheorem{lemma}[definition]{Lemma}

\newenvironment{prof7}[1][Proof]{\textbf{#1.} }{\ \rule{0.5em}{0.5em}}

\begin{document}
\title[Compactness of isospectral conformal Finslerian metrics set on a $3$-manifold]{Compactness of isospectral conformal Finslerian metrics set on a $3$-manifold} % please provide
                                % an abbreviated title 
\author{G. Nibaruta and P. Nshimirimana}

\address{Ecole Normale Sup\'erieure\\ 
\noindent
P.~O. Box 6983 \\ 
 Bujumbura\\ Burundi}

 \email{nibarutag@gmail.com}

\begin{abstract}
Let $F$ be a Finslerian metric on an $n-$dimensional closed manifold $M$. In this work, we study problems about compactness 
 of isospectral sets of  conformal Finslerian  metrics when $n=3$.% More precisely, let 
 %$\widetilde{F}\in[F]$ be a Finslerian metric in the conformal class of $F$ whose scalar curvature is nonpositive constant. We show that the 
%set of metrics in $[F]$  isospectral to $\widetilde{F}$ is compact in the $C^{\infty}$-topology.
\end{abstract} 

\keywords{Finslerian conformal  metrics; isospectral metrics; %\\
%\hspace*{2.3cm}~~~~sphere bundle, Finsler-Ehresmann connection, 
Sobolev spaces for manifolds.}

\maketitle

\section{Introduction}\label{Section1}
The spectral theory of a  Laplacian associated with a metric is a very important subject
in Differential Geometry \cite{bi2017}. For example, to know the Laplacian spectrum of a Riemannian metric $g$ makes it possible to know the geometry of the Riemannian manifold. However, there are only a few examples of manifolds whose the 
spectrum is known explicitly. More so, spectral problems
have not yet been well studied for Finslerian  metrics even if these last are in these years very important.

Consider a $C^{\infty}$ manifold $M$. Let $F$ and $\widetilde{F}$ be two Finslerian metrics on $M$. $F$ and $\widetilde{F}$ are said to be \textit{conformal} if there exists a $C^{\infty}$ function $u$ on $M$ such that $\widetilde{F}=e^{u}F$. Then, $F$ and $\widetilde{F}$ 
 are said to be conformal. 
 Two Finslerian metrics $F$ and $\widetilde{F}$ on a compact manifold are said to be
\textit{isospectral} if their associated \textit{Laplacian operators}, $\hat{\Delta}$ and $\hat{\widetilde{\Delta}}$, 
on functions, have the same spectrum. If, in addition, $F$ and $\widetilde{F}$ are conformal to each
other, they will be said to be isospectral conformal metrics.

We are going to answer to the following questions:
Do exist Finslerian metrics $F$ and $\widetilde{F}$ on $M$ whose $F$ and $\widetilde{F}$ 
 are conformally related non-isometric such that $Spec(\hat{\widetilde{\Delta}})=Spec(\hat{\Delta})$ ?
 If the answer question is yes, are the conformal factors bounded ?

To answer to these questions, we consider a $3-$dimensional closed Finslerian manifold $(M,F)$ and its sphere bundle $SM:=\cup_{x\in M}S_xM$
 where $S_xM$ is the set consisting of all directions $[y]:=\{ty:t>0, y\in T_xM\backslash\{0\}\}$. From \cite{bi1}, $SM$ is 
a $5-$dimensional closed Riemannian manifold of all constant direction tangent on $3-$dimensional closed Finsler manifold $M$.
Then the functions defined on Finsler manifold $(M,F)$ will be seen as the functions defined on $SM$ with constant direction in $S_xM$.
%because 
%the conformal factors induced by the conformal deformation of $F$ do not depend on the direction as shown by the following result:
%\begin{proposition}
 %Let $F$ and $\widetilde{F}$ be two Finsler metrics on a smooth manifold $M$. If $F$ and $\widetilde{F}$ are conformal, 
 %then the conformal factor $\varphi\in C^\infty(\mathring{T}M)$ is not depend of direction.
%\end{proposition}

%  It is known \cite{bi2}, for a Riemannian manifold, when $(M^3,g)$ has constant negative scalar curvature, an isospectral set of
% Riemannian metrics $\widetilde{g}={\varphi}^4g$ conformal to $g$ is compact in the $C^\infty$ topology.
 %In the second section, we recall some basic notions which are necessarly in the following. 
 %In section 3, we give the definition of the horizontal Laplacian $\hat{\Delta}$ operator, defined in \cite{bi5}
%   for $C^{\infty}$ functions defines on Finsler manifolds. 
%   We propose some properties on $\hat{\Delta}$ in order to show the existence of spectrum of it.
%  In section 4, motivated by Sunada method, we construct Finsler non-isometric manifolds. As application, using the Yamabe (type) equation
%   found for Finsler manifolds $(M,F)$ when the Finsler-Ehresmann connection $\mathcal{H}$ is invariant under the conformal change 
%   $\widetilde{F}=\varphi^{\frac{2}{n-2}}F$ ($\varphi\in C^{\infty}(SM)$, $n=3$), we prove, in section 5,
%   that the set of Finsler metrics in the conformal
% class $[F]$ that are isospectral to $\widetilde{F}$ is compact in the $C^2$-topology.

 The outline of this document is organized as follows. 
  In Section \ref{Section2},  
  we establish some basic notions on Finslerian Spectral Geometry.
  The Section \ref{Section3} is devoted to conformal deformations of some  quantities related with Finslerian metrics. %In the Section \ref{Section4} we construct isospectral Finsler manifolds non-isometric.
  Then we prove in the Section \ref{Section5} the main results.\\
  
  \noindent
\textbf{Notation.}
 In the following, $M$ is an $3-$dimensional manifold of class $C^2$. We denote by $T_xM$ the tangent space at 
   $x\in M$ and by $TM:=\bigcup_{x\in M}T_xM$ the tangent bundle of $M$.
   Set $\mathring{T}M=\{(x,y)\in TM: y\neq 0\text{ where } 0 \text{ is a zero vector } \in T_{x}M\}$ and $\pi:TM\longrightarrow M:\pi(x,y)\longmapsto x$ the natural projection. 
   Let $(x^1,...,x^n)$ be a local coordinate 
   on an open subset $U$ of $M$ and $(x^1,...,x^n,y^1,...,y^n)$ be the local coordinate 
   on $\pi^{-1}(U)\subset TM$, the local coordinate system $(x^i)_{i=1,...,n}$ produces the coordinate bases  
   $\{\frac{\partial}{\partial x^i}\}_{i=1,...,n}$ and $\{dx^i\}_{i=1,...,n}$ respectively, for 
   $TM$ and cotangent bundle $T^*M$.\\
  
  \noindent
\textbf{Convention.}
   We use Einstein summation convention: repeated
   upper and lower indices will automatically be summed, from $1$ to $n$, unless otherwise will be noted. For example 
   $g_{ij}v^iv^j=\sum_{i,j=1}^n g_{ij}v^iv^j.$
  
 \section{Some basic notions on Finslerian Spectral Geometry}\label{Section2} % The names of sections should be written in sentence-case

   \begin{definition}\label{defi1}Let $M$ be an $n$-dimensional manifold. A function $F:TM\longrightarrow [0,\infty)$ is a \textit{Finslerian metric} 
   	on $M$ if :
   	\begin{itemize}
   		\item [(i)] $F$ is $C^{\infty}$ on the entire slit tangent bundle $\mathring{T}M$,
   		\item [(ii)] $F$ is positively $1$-homogeneous on the fibers of $TM$, that is \\
   		$\forall c>0,~F(x,cy)=cF(x,y),$
   		\item [(iii)] the Hessian matrix $(g_{ij}(x,y))_{1\leq i,j\leq n}$ with elements
   		\begin{eqnarray}\label{1}
   		g_{ij}(x,y):=\frac{1}{2}\frac{\partial^2F^2(x,y)}{\partial y^i\partial y^j}
   		\end{eqnarray}
   		is positive-definite at every point $(x,y)$ of $\mathring{T}M$.
   	\end{itemize}
   \end{definition}
\begin{remark}
	The Hessian matrix $g$ whose elements are defined in (\ref{defi1}) is a natural Riemannian metric on the pulled-back bundle $\pi^*TM$ over the manifold 
	$\mathring{T}M$. $\pi^*TM$ is a vector bundle and its fiber at a point $(x,y)\in\mathring{T}M$ is
\begin{eqnarray}\label{12}
 \pi^*TM |_{(x,y)}:=\{(x,y;v): v\in T_xM\}|_{(x,y)}\cong T_xM.
\end{eqnarray}
\end{remark}
\subsection{Non-linear  connection on the slit tangent bundle}
Consider the differential map $\pi_*$  of the submersion 
$\pi:\mathring{T}M \longrightarrow M$.
The vertical subspace of $T\mathring{T}M$ is defined by 
$\mathcal{V}:=ker({\pi_*})$ 
and is locally spanned by the set $\{F\frac{\partial}{\partial y^i}, 1\leq i\leq n\}$,
on each $\pi^{-1}(U)\subset \mathring{T}M$.

An horizontal subspace $\mathcal{H}$ of $T\mathring{T}M$ is by definition any complementary to 
$\mathcal{V}$. The bundles $\mathcal{H}$ and $\mathcal{V}$ give a smooth splitting 
\begin{eqnarray}\label{decomposition}
T\mathring{T}M=\mathcal{H}\oplus\mathcal{V}.  
\end{eqnarray}
An non-linear connection is a selection of a horizontal subspace $\mathcal{H}$ of $T\mathring{T}M$.
%It is known \cite{Chern} that $\mathcal{H}$ can be canonically defined from the geodesics equation.

\begin{definition}\label{defi2} Let $\pi:\mathring{T}M\longrightarrow M$ be the 
	restricted projection.
	\begin{itemize}
		\item [(1)]A non-linear connection of $\pi$
		is the subbundle $\mathcal{H}$ of $T\mathring{T}M$ given by 
		\begin{eqnarray}
		\mathcal{H}:=ker \theta,
		\end{eqnarray}
		where $\theta:T\mathring{T}M\longrightarrow \pi^*TM$ is the bundle morphism defined 
		by 
		\begin{eqnarray}
		\label{03b}
		\theta|_{(x,y)}=\frac{\partial}{\partial x^i}\otimes \frac{1}{F}(dy^i+N_j^idx^j)
		\end{eqnarray}
		with $N_j^i(x,y):=\frac{\partial G^i(x,y)}{\partial y^j}$ 
		for 
		\begin{eqnarray}\label{Gcoef}
		G^i(x,y):=\frac{1}{4}g^{il}(x,y)
		\Big[\frac{\partial g_{jl}}{\partial x^k}(x,y)+\frac{\partial g_{kl}}{\partial x^j}(x,y)-\frac{\partial g_{jk}}{\partial x^l}(x,y)\Big]y^jy^k.
		\end{eqnarray}
		\item [(2)]The form $\theta:T\mathring{T}M\longrightarrow \pi^*TM$ induces a linear map 
		\begin{eqnarray}
		\theta|_{(x,y)}:T_{(x,y)}\mathring{T}M\longrightarrow T_xM,
		\end{eqnarray} 
		for each point $(x,y)\in \mathring{T}M$; where $x=\pi(x,y)$. \\
		The vertical lift of a section $\xi$ of $\pi^*TM$ is a unique section $\textbf{v}(\xi)$ of $T\mathring{T}M$ 
		such that for every $(x,y)\in \mathring{T}M$, 
		\begin{eqnarray}
		\pi_* (\textbf{v}(\xi))|_{(x,y)}=0_{(x,y)}\text{     and     }\theta (\textbf{v}(\xi))|_{(x,y)}=\xi_{(x,y)}.
		\end{eqnarray}
		\item [(3)]The differential projection $\pi_*:T\mathring{T}M\longrightarrow \pi^*TM$ induces a linear map 
		\begin{eqnarray}
		\pi_*|_{(x,y)}:T_{(x,y)}\mathring{T}M\longrightarrow T_xM,
		\end{eqnarray} 
		for each point $(x,y)\in \mathring{T}M$; where $x=\pi(x,y)$. \\
		The horizontal lift of a section $\xi$ of $\pi^*TM$ is a unique section $\textbf{h}(\xi)$ of $T\mathring{T}M$ 
		such that for every $(x,y)\in \mathring{T}M$,
		\begin{eqnarray}
		\pi_* (\textbf{h}(\xi))|_{(x,y)}=\xi_{(x,y)}\text{     and     }\theta (\textbf{h}(\xi))|_{(x,y)}=0_{(x,y)}.
		\end{eqnarray}
	\end{itemize}
\end{definition}

We have the following.
\begin{definition}\label{defi3b} 
	A Finslerian tensor field $T$ of type $(q,0;p_1,p_2)$
	on $\mathring{T}M$ is a $C^{\infty}$ section of the tensor bundle 
	\begin{eqnarray}
	\underbrace{\pi^*T^*M\otimes...\otimes\pi^*T^*M}_{p_1-times}
	\otimes \underbrace{T^*\mathring{T}M\otimes...\otimes T^*\mathring{T}M}_{p_2-times}
	\otimes\bigotimes^q\pi^*TM.  
	\end{eqnarray}
\end{definition}

\begin{remark}\label{rm2}
	In a local chart, $$T=T_{i_1...i_{p_1}j_1...j_{p_2}}^{k_1...k_q}\partial_{k_1}\otimes ... \otimes\partial_{k_q}
	\otimes dx^{i_1}\otimes...\otimes dx^{i_{p_1}}\otimes \varepsilon^{j_{1}}\otimes ... \otimes\varepsilon^{j_{p_2}}$$
	where $(\partial_{k_1}\otimes ... \otimes\partial_{k_q}
	\otimes dx^{i_1}\otimes...\otimes dx^{i_{p_1}}\otimes \varepsilon^{j_{1}}\otimes ... \otimes\varepsilon^{j_{p_2}})
	_{k\in\{1,...,n\}^q, i\in\{1,...,n\}^{p_1}, j\in\{1,...,n\}^{p_2}}$ is a basis section of this tensor 
	and, the $\partial_{k_r}:=\frac{\partial}{\partial x^{k_r}}$ as well as $\varepsilon^{j_{s}}$
	are respectively the basis sections for $\pi^*TM$ and $T^*\mathring{T}M$ dual of $T\mathring{T}M$. 
\end{remark}
\begin{example}
	\begin{itemize}
		\item [(1)]The Hessian matrix $g$, defined in (\ref{1}),  is of type $(0,0;2,0)$.
		\item [(2)] The Ehresmann-Finsler form $\theta$ is of type $(1,0;0,1)$.
	\end{itemize}
\end{example}
\subsection{Some linear connection on the vector bundle $\pi^*TM$}
The following lemma defines the Chern connection on $\pi^*TM$.
\begin{lemma}\cite{bi5}\label{lem1}
	Let $(M,F)$ be a Finslerian manifold and $g$ its fundamental tensor. 
	There exists a unique linear connection $\nabla$ 
	on the vector bundle $\pi^*TM$ such that, for all 
	$X,Y\in \chi(\mathring{T}M)$ and for every $\xi, \eta\in\Gamma(\pi^*TM)$, one has the following properties:
	\begin{itemize}
		\item [(i)]
		$\nabla_X\pi_*Y-\nabla_Y\pi_*X=\pi_*[X,Y],$
		\item [(ii)] 
		$X(g(\xi,\eta))=g(\nabla_X\xi,\eta)+g(\xi,\nabla_X\eta)+2\mathcal{A}(\theta(X),\xi,\eta)$ \\
		where $\mathcal{A}:=\frac{F}{2}\frac{\partial g_{ij}}{\partial y^k}dx^i\otimes dx^j\otimes dx^k$ is the Cartan tensor.
	\end{itemize}
\end{lemma}

One has $\nabla_{\frac{\delta}{\delta x^j}}\frac{\partial}{\partial x^k}=\Gamma_{jk}^i\frac{\partial}{\partial x^i}$ where
\begin{eqnarray}\label{Ccoef}\Gamma_{jk}^i&:=&
\frac{\partial^2 G^i}{\partial y^j\partial y^k}
\end{eqnarray}
which can be written as
\begin{eqnarray}
\Gamma_{jk}^i&=&\frac{1}{2}g^{il}\left(\frac{\delta g_{jl}}{\delta x^k}+\frac{\delta g_{lk}}{\delta x^j}
-\frac{\delta g_{jk}}{\delta x^l}\right)\label{8775244}
\end{eqnarray}
with
\begin{eqnarray}
\left\{\frac{\delta}{\delta x^i}:=\frac{\partial}{\partial x^i}-N_i^j\frac{\partial}{\partial y^j}
=\textbf{h}(\frac{\partial}{\partial x^i})\right\}_{i=1,...,n}.\label{562564654}
\end{eqnarray}
\subsection{Curvatures}
\begin{definition}
	The full curvature associated with the Chern connection $\nabla$
	on the vector bundle $\pi^*TM$ over the manifold $\mathring{T}M$ is the application
	$$ \phi:\begin{matrix}
	\chi(\mathring{T}M)\times\chi(\mathring{T}M)\times\Gamma(\pi^*TM)&\to&\Gamma(\pi^*TM)\\
	(X,Y,\xi)&\mapsto
	&\phi(X,Y)\xi=\nabla_X\nabla_Y\xi-\nabla_Y\nabla_X\xi-\nabla_{[X,Y]}\xi.
	\end{matrix}$$
\end{definition}
By the relation (\ref{decomposition}), we have
\begin{eqnarray}
\nabla_X=\nabla_{\hat{X}}+\nabla_{\check{X}},
\end{eqnarray}
where $X=\hat{X}+\check{X}$ with $\hat{X}\in\Gamma(\mathcal{H})$ and $\check{X}\in\Gamma(\mathcal{V})$.

Using the metric $F$, one can define the full curvature  
of $\nabla$ as:
\begin{eqnarray}
\Phi(\xi,\eta,X,Y)&=&g(\phi(X,Y)\xi,\eta)\nonumber\\
&=&g(\phi(\hat{X},\hat{Y})\xi+\phi(\hat{X},\check{Y})\xi
+\phi(\check{X},\hat{Y})\xi+\phi(\check{X},\check{Y})\xi,\eta)\nonumber\\
&=&\textbf{R}(\xi,\eta,X,Y)+\textbf{P}(\xi,\eta,X,Y)+\textbf{Q}(\xi,\eta,X,Y),
\end{eqnarray}
where 
$
\textbf{R}(\xi,\eta,X,Y)=g(\phi(\hat{X},\hat{Y})\xi,\eta),$ 
$\textbf{P}(\xi,\eta,X,Y)=g(\phi(\hat{X},\check{Y})\xi,\eta)
+g(\phi(\check{X},\hat{Y})\xi,\eta)$ and $\textbf{Q}(\xi,\eta,X,Y)=g(\phi(\check{X},\check{Y})\xi,\eta)$
are respectively the \textit{first (horizontal) curvature}, 
\textit{mixed curvature} and \textit{vertical} curvature.

In particular, if $\nabla$ is the Chern connection, the $\textbf{Q}$-curvature vanishes.

\begin{definition}
 Let $F$ be a Finslerian metric, $\textbf{R}$ the first curvature associated with $F$ and $\{e_a\}_{a=1,...,n}$ an $g$-orthonormal basis sections of $\pi^*TM$. 
  \begin{itemize}
   \item [(1)] The first Ricci curvature $\hat{\textbf{r}}$ of $(M,F)$ is defined by 
   \begin{eqnarray}
    \hat{\textbf{r}}(\xi,X):=trace_g(\eta\longmapsto \textbf{R}(X,\textbf{h}(\eta)\xi), ~~\forall X\in\Gamma(T\mathring{T}M)~~and~~\xi,\eta\in\Gamma(\pi^*TM).
   \end{eqnarray}
In $g$-orthonormal basis sections $\{e_a\}_{a=1,...,n}$ of $\pi^*TM$, one has 
\begin{eqnarray}
    \hat{\textbf{r}}(\xi,X):=\sum_{a=1}^n\textbf{R}(e_a,\xi,X,\textbf{h}(e_a)).
   \end{eqnarray}
   \item [(2)] The first  scalar curvature $\hat{\textbf{S}}$ associated with $F$ is the trace of the first  Ricci curvature. 
   $\hat{\textbf{S}}$ is a function on $\mathring{T}M$ or on $SM$. In $g$-orthonormal basis sections $\{e_a\}_{a=1,...,n}$ of $\pi^*TM$, one has 
\begin{eqnarray}
    \hat{\textbf{S}}:=\sum_{a=1}^n\hat{\textbf{r}}(e_a,\textbf{h}(e_a))=\sum_{a,b=1}^n\textbf{R}(e_a,e_b,\textbf{h}(e_b),\textbf{h}(e_a)).
   \end{eqnarray}
  \end{itemize}
\end{definition}
\subsection{The Finslerian covariante derivative% by the pulled-back bundle approach
}
\begin{definition}%\textbf{$($Chern horizontal covariant derivative$)$}.
 Let $(M,F)$ be a Finslerian manifold %, $SM$ its sphere bundle
 and $\nabla$ 
 the Chern connection on pulled-back tangent bundle $\pi^*TM$. We denote by $\hat{\nabla}$ 
 the Chern connection in the horizontal direction given by:
 $$\hat{\nabla}\xi(\eta):=\nabla_{\textbf{h}(\eta)}\xi,~~~\forall \xi, \eta\in\Gamma(\pi^*TM).$$
 For a $C^\infty$-function $u$ on $M$, one defines the Chern 
 horizontal covariant derivative $\hat{\nabla} u$ of $u$ by the following expression
 $$(\hat{\nabla}u)(\xi):=\nabla_{\textbf{h}(\xi)}u:=du(\textbf{h}(\xi)), \xi\in\Gamma(\pi^*TM).$$
 \end{definition}
 \begin{remark}\label{0012}
  From the fact that all objects involved are of constant horizontal direction, we can see $u\in C^\infty(M)$
   as $u\in C^\infty(SM)$ of constant direction $[y]$.
 \end{remark}
 Then for $u\in C^{\infty}(SM)$ second Chern horizontal covariant derivative is define by 
 \begin{eqnarray}
  (\hat{\nabla}^2u)(\xi,\eta)&:=&(\nabla_{\textbf{h}(\xi)}\hat{\nabla}u)(\eta),\nonumber\\
                           %&=&\nabla_{\textbf{h}(\xi)}(\hat{\nabla}u(\eta))-(\hat{\nabla}u)(\nabla_{\textbf{h}(\xi)}\eta),\nonumber\\
                           &=&\nabla_{\textbf{h}(\xi)}\nabla_{\textbf{h}(\eta)}u-\hat{\nabla}u(\nabla_{\textbf{h}(\xi)}\eta).
 \end{eqnarray}
By extension, for all $k\in\mathbb{N}$, the $k-$th Chern covariant derivative in the horizontal direction of $u$
is denoted by 
$$\hat{\nabla}^ku:=\hat{\nabla}(\hat{\nabla}^{k-1}u).$$
In particular, $\hat{\nabla}^0u=u$. One denotes by $||\hat{\nabla}^ku||$ the norm of $\hat{\nabla}^ku$ 
induced by the fundamental tensor $g$ of $F$. It is given by 
$$||\hat{\nabla}^ku||^2=g^{i_1j_1} ... g^{i_kj_k}(\hat{\nabla}^ku)_{i_1 ... i_k}(\hat{\nabla}^ku)_{j_1 ... j_k}.$$

\subsection{The Finslerian Laplacian% by the pulled-back bundle approach
}
\begin{definition}\label{defi2}
	Let $(M,F)$ be an $n$-dimensional Finslerian manifold and $g$ the Hessian matrix associated with $F$.
	The gradient of a smooth function $u$ on $M$ is, the section of the vector bundle $\pi^*TM$ denoted by  $\triangledown u$, 
	 given by
	 \begin{eqnarray}
	 	g_{(x,y)}(\triangledown u_{(x,y)},\xi_{(x,y)})=du_{(x,y)}(\xi_{(x,y)})
	 \end{eqnarray}
	 for any $\xi\in\Gamma(\pi^*TM)$ and for every $(x,y)\in \mathring{T}M$. Locally, 
	 \begin{eqnarray}
	 \triangledown u_{(x,y)}=g^{ij}(x,y)\frac{\partial u}{\partial x^i}\frac{\partial }{\partial x^j}.
	 \end{eqnarray}
\end{definition}
 \begin{definition}%\textbf{$($Horizontal divergence$)$}\cite{bi5}.
Let $(M,F)$ be a Finsler manifold. One define the horizontal divergence of a vector field $\xi\in\Gamma(\pi^*TM)$ by 
$$\hat{D}\xi=trace_g(\eta\longmapsto\nabla_{\textbf{h}(\eta)}\xi)$$
% and the vertical divergence of a vector field $\xi\in\Gamma(\pi^*TM)$
% $$div^v\xi=trace_g(\eta\longmapsto\nabla_{\eta^v}\xi)$$
where $g$ is a fundamental tensor and $\nabla$ is the Chern connection. 
 \end{definition}
 \begin{remark}
  In the basis sections $\{\frac{\partial}{\partial x^i}\}_{i=1,...,n}$ of the pulled-back tangent bundle $\pi^*TM$, one get 
\begin{eqnarray}
\hat{D}\xi=g^{ij}g\left(\nabla_{\frac{\delta}{\delta x^i}}\xi,\frac{\partial}{\partial x^j}\right).
\end{eqnarray}
 \end{remark}
 
Now, one defines a Hessian of the $C^\infty$ function $u$, denoted by $Hu$, 
  on $(M,F)$ as
  \begin{definition}%\textbf{$($Hessien$)$}\cite{biMbatakou2014}.
  Suppose $u\in C^{\infty}(M)$. A Hessian of $u$ is the mapping
  \begin{eqnarray}
   Hu:\Gamma(\pi^*TM)\times \Gamma(T\mathring{T}M)\longrightarrow C^{\infty}(\mathring{T}M): 
   Hu(\xi,X)=g\left(\xi, \nabla_X(\triangledown u)\right).
  \end{eqnarray}
 \end{definition}
\begin{definition}%\textbf{$($Horizontal Laplacian, Vertical Laplacian$)$}\cite{biMbatakou2014}.
Let $(M,F)$ be a $C^{\infty}$ Finsler manifold and $u$ a $C^{\infty}$ function on $M$.
% With identification $\pi^*TM\cong\mathcal{H}$, 

The horizontal Laplacian 
$\hat{\Delta}u$ of $u$ is given by 
\begin{eqnarray}\label{lp}
\hat{\Delta}u=-\hat{D}\triangledown u
\end{eqnarray}
and %with identification $\pi^*TM\cong\mathcal{V}$, 
the vertical Laplacian $\check{\Delta}u$ of $u$ by 
\begin{eqnarray}
 \check{\Delta}u=-\check{D}\triangledown u.
\end{eqnarray}
 \end{definition}
\begin{proposition}\cite{bi5}
   Let $u\in C^{\infty}(M)$. The horizontal Laplacian $\hat{\Delta}u$ and the vertical Laplacian $\check{\Delta}u$ of
   $u$ can be given in term of the Hessian of $u$ by
    \begin{eqnarray}
    \label{0}
    \hat{\Delta} u=-trace_g\left((\xi,\eta)\mapsto Hu(\xi,\textbf{h}(\eta))\right),~~~~  \xi,\eta\in\Gamma(\pi^*TM),
    \end{eqnarray}
    and 
    \begin{eqnarray}
    \check{\Delta} u=-trace_g\left((\xi,\eta)\mapsto Hu(\xi,\textbf{v}(\eta))\right),~~~~  \xi,\eta\in\Gamma(\pi^*TM).
    \end{eqnarray}
   Furthermore, in $g-$orthonormale basis sections $\{e_a\}_{a=1,...,n},$  one has
    \begin{eqnarray}
    \label{000}
     \hat{\Delta} u=-\sum_{a=1}^nHu(e_a,\textbf{h}(e_a))=-\hat{\nabla}^a\hat{\nabla}_au.
    \end{eqnarray}
    and 
    \begin{eqnarray}
     \check{\Delta} u=-\sum_{a=1}^nHu(e_a,\textbf{v}(e_a)).
    \end{eqnarray}
  \end{proposition}
  
 \section{Conformal deformations of some  quantities related with Finslerian metrics}\label{Section3}
       
       In this paper, we do not recall the notions about tensor formalism in Finsler geometry. They have been studied in \cite{bi5}. 
%  There is no unified definition of scalar curvature in Finsler geometry, although 
%  several geometers have offered several versions of the definition of scalar curvature.
We adopt the definition of horizontal scalar curvature introduced in \cite{bi21} and \cite{bi5}.
 %We consider, a horizontal scalar curvature $Scal^h$ as function on the sphere bundle $SM$ over a Finsler manifold $(M^n,F)$.
%Set $\hat{\textbf{S}}=Scal^h$. We have, in $g-$orthonormal basis sections $\{e_a\}_{a=1,...,n}$ for $\pi^*TM$,
 %\begin{eqnarray}
 % \hat{\textbf{S}}:=\sum_{a=1}^nRic^h(e_a,e_a^h):=\sum_{a,b=1}^nR(e_b,e_a,e_a^h,\textbf{h}(e_b))
 %\end{eqnarray}
%where $Ric^h$ and $R$ are, respectively, horizontal Ricci tensor and horizontal part of the full curvature tensor associated to the 
%the Chern connection.

\subsection{Conformal change of volume element and horizontal scalar curvature}
Let $(M^n,F)$ be a compact $C^{\infty}$ Finsler manifold of dimension $n\geq3$, $SM$ its sphere bundle and 
 $\hat{\textbf{S}}$ the first scalar curvature associated to $F$. Given a $C^{\infty}$ function $u:SM\longrightarrow \mathbb{R}$ 
 on $SM$  with $u(x,[y])\neq0$ for all $(x,[y])\in SM$ let $\widetilde{F}=e^uF$ be the conformal change of the Finsler metric $F$.
 
 \subsubsection{Conformal change of volume}
 Under the conformal change $\widetilde{F}=e^uF$ of the form $\widetilde{F}=\varphi^{\frac{2}{n-2}}F$, where $\varphi$ 
 is $C^{\infty}$ function on $SM$ with $\varphi(x,[y])\neq0 $, 
 the fundamental tensors $g$ and $\widetilde{g}$ associated the respective Finsler metrics 
$F$ and $\widetilde{F}$ are conformally related by 
\begin{eqnarray}
 \widetilde{g}=\varphi^{\frac{4}{n-2}}g.
\end{eqnarray}
while the volume forms
 on $SM$ associated to the Finsler metrics $\widetilde{F}$ and $F$ are related by 
 \begin{eqnarray}
  \eta_{\widetilde{F}}=\varphi^{\frac{2n}{n-2}}\eta_F.
 \end{eqnarray}
For a $3-$dimensional manifold $M$, we have 
$$\widetilde{F}=\varphi^{2}F, ~~\widetilde{g}=\varphi^{4}g, ~~\eta_{\widetilde{F}}=\varphi^{6}\eta_F.$$
\subsubsection{Conformal change of horizontal scalar curvature}
 From \cite{bi5}, if $\widetilde{F}=e^uF$ is a conformal metric of $F$ then the associated 
 horizontal scalar curvature $\hat{\widetilde{\textbf{S}}}$ will be related to $\hat{\textbf{S}}$ by the equation
 \begin{eqnarray}
  \hat{\widetilde{\textbf{S}}} &=& e^{-2u}\left[\hat{\textbf{S}}+2(n-1)\hat{\Delta}u-(n-1)(n-2)||\triangledown u||_g^2\right]\nonumber\\
                  & &+e^{-2u}\left[\sum_{b=1}^{n}(2-n)[du(\Theta(\textbf{h}(e_b),\textbf{h}(e_b)))+g(\Theta(\textbf{h}(e_b),\hat{\nabla}u),\textbf{h}(e_b))]\right]\nonumber\\
 & &+e^{-2u}\left[\sum_{a,b=1}^{n}[g(\Theta(\textbf{h}(e_b),\Theta^h(\textbf{h}(e_b),\textbf{h}(e_b))),e_a)-g(\Theta(e_a^h,\Theta^h(\textbf{h}(e_b),\textbf{h}(e_b))),e_a)]\right]\nonumber\\
& &+e^{-2u}\left[\sum_{a,b=1}^{n}[g((\nabla_{\textbf{h}(e_b)}\Theta)(\textbf{h}(e_b),\textbf{h}(e_b)),e_a)-g((\nabla_{e_a^h}\Theta)(\textbf{h}(e_b),\textbf{h}(e_b)),e_a)]\right].
 \end{eqnarray}
 Where $\Theta$ is a $(0,2;1)-tensor$ defined in \cite{bi5}, $\{e_a\}_{a=1,...,n}$ is a $g-$orthonormale basis sections for 
 $\pi^*TM$, defined in  \cite{bi1}.
Using the results in \cite{bi5}, when the Finsler-Ehresmann connection is invariant under the conformal change of 
metric then 
\begin{eqnarray}
 \hat{\widetilde{\textbf{S}}} &=& e^{-2u}[\hat{\textbf{S}}+2(n-1)\hat{\Delta}u-(n-1)(n-2)||\triangledown u||_g^2].
\end{eqnarray}

If we consider the conformal deformation in the form $\widetilde{F}=\varphi^{\frac{2}{n-2}}F$ (with $\varphi\in C^{\infty},
\varphi>0)$, the horizontal scalar curvature $\hat{\widetilde{\textbf{S}}}$ associated to $\widetilde{F}$ satisfies the equation
\begin{eqnarray}
\label{11}
 \frac{4(n-1)}{n-2}\hat{\Delta}\varphi+\hat{\textbf{S}}\varphi=\hat{\widetilde{\textbf{S}}}\varphi^{\frac{n+2}{n-2}}
\end{eqnarray}
where $\hat{\Delta}\varphi:=-\hat{D}\triangledown\varphi$ is the horizontal Laplacian of $\varphi$.
When $M$ is of three dimensional, we have the Yamabe type equation
\begin{eqnarray}
\label{22}
 8\hat{\Delta}\varphi+\hat{\textbf{S}}\varphi=\hat{\widetilde{\textbf{S}}}\varphi^5.
\end{eqnarray}
\begin{remark}
 According to \cite{bi11}, we have the equation
 \begin{eqnarray}
\label{2}
 \frac{4(N-1)}{N-3}\hat{\Delta}\varphi-\hat{\textbf{S}}\varphi=-\hat{\widetilde{\textbf{S}}}\varphi^{\frac{N+5}{N-3}}
\end{eqnarray}
when $N=2n-1$ is the dimension of the Riemannian manifold $SM$.
\end{remark}

\subsection{Heat kernel asymptotics of horizontal Laplacian $\hat{\Delta}$}
Suppose an $n-$dimensional compact Finsler manifold $(M,F)$. Among the main tools used in studying the spectral compactness 
is the heat trace expansion.\\

\noindent
\subsubsection{For Riemannian manifolds} When $(M,F)$ is Riemannian, that is $F=\sqrt{g_x}$, it is 
well known that \cite{bi2017}
% Let $H(t,x,z)$ be the fundamental solution to the heat equation (or the heat kernel)
%  $$\frac{\partial H}{\partial t}+\Delta_{g_x}H=0,~~~H:\mathbb{R}^+\times M\longrightarrow\mathbb{R},$$
%  where $(t,x)\in \mathbb{R}^+\times M$, $\Delta_{g_x}$ is the Laplace-Beltrami operator on Riemannian manifold $(M,g_x)$.
%  It is known, from \cite{bi2}, the asymptotic developpement for the heat kernel associated to the 
%  Laplace-Beltrami operator $\Delta_{g_x}$ with respect to the Riemannian metric $g_x$ is
%  If we denote by $0<\lambda_0<\lambda_1<\lambda_2,...$ the eigenvalues of the Laplacian $\Delta_{g_x}$, 
% then, it is known  the trace of $H(t,x,z)$ has the following asymptotic expansion:
\begin{eqnarray}
%  Tr(H(t,x,x))&=&\int_MH(t,x,x)dv_{g_x}\nonumber\\
Tr(e^{-t\Delta_{g_x}}) &\sim&\frac{1}{(4\pi t)^{\frac{n}{2}}}(a_0+a_1t+a_2t^2+ ...),
\end{eqnarray}
as $t\rightarrow 0$, where $dv_{g_x}$ denote the Riemannian volume element of $(M,g_x)$ 
and $\Delta_{g_x}$ is the Laplace-Beltrami operator on $(M,g_x)$.
 The heat invariants $a_0, a_1, a_2, ...$ are integrals of derivatives of curvature terms on Riemannian 
manifold $M$.\\

\noindent
\subsubsection{For Finsler manifolds}
 Recall that, the sphere bundle $SM$ over $(M,F)$ is compact Riemannian manifold %of dimension $N=(2n-1)$
\cite{bi1}.
We may use a conformal factor $\varphi$ and the horizontal quantities given above to express 
the first several heat invariants, for 3-dimensional closed Finsler manifold, as integrals of horizontal 
curvature terms on Riemannian manifold $SM$. They are explicitely given by 
\begin{eqnarray}
 \widetilde{a_0}&=&\int_{SM}\eta_{\widetilde F}\nonumber\\
    %&=&\int_{SM}\varphi^{\frac{2N+2}{N-3}}\eta_ F\nonumber\\
    &=&\int_{SM}\varphi^{6}\eta_ F,
\end{eqnarray}
\begin{eqnarray}
 \widetilde{a_1}&=&\frac{1}{6}\int_{SM}\hat{\widetilde{\textbf{S}}}\eta_{\widetilde F}\nonumber\\
    &=&\frac{1}{6}\int_{SM}(\hat{\textbf{S}}\varphi^2+8|\hat{\nabla} \varphi|^2)\eta_ F,
\end{eqnarray}
\begin{eqnarray}
 \widetilde{a_2} =\frac{1}{360}\int_{SM}(3(\hat{\widetilde{\textbf{S}}})^2+6|\widetilde{\rho}^h|^2)\eta_{\widetilde F},\\
\end{eqnarray}
when $dim SM=5$ where $|\widetilde{\rho}^h|=||\widetilde{R}^h||,
\widetilde{R}^h$ the horizontal Ricci tensor associated to the Finsler metric 
$\widetilde{F}$.

\section{Main results}\label{Section5}
\subsection{Non-isometric isospectral Finsler manifolds}\label{Section4}
%In this section we use Sunada method \cite{bi10} to construct isospectral Finsler manifolds non-isometric. 
%Let consider $(M,F)$ a 3-Finsler manifold and $SM$ its sphere bundle of constant directions.
% The Sunada method reduce the construction of isospectral manifolds to a finite group problem. 
One has the following result
\begin{theorem}
Let $(M,F)$ be a connected closed 3-Finsler manifold and $SM$ its sphere bundle of constant directions $y$. Then there exist, 
non-isometric,
conformal Riemannian metrics $g^s$ and $\widetilde{g}^s\in[g^s]$, on $SM$ such that $(SM,g^s)$ and 
$(SM,\widetilde{g}^s)$ are isospectral.
\end{theorem}

\begin{prof7}
Since $g^s$ is a Riemannian metric, the proof is similar to the proof found in \cite{bi2}.
\end{prof7}

In what follows, all Finsler manifolds $M$ will be closed and of 3-dimensional.

\subsection{Sobolev spaces on the space $SM$} 
% \begin{definition}
 Let $(M,F)$ be a smooth Finsler manifold of dimension $n$ and $SM$ the space of all constant directions tangent on $M$.
 Given a $C^{k}$ function $\varphi$ on $SM (k\geq0$ an integer $)$, we denote by $\hat{\nabla}^k$ 
 the $k$-th horizontal covariant derivative of $\varphi$ define in paragraph $2.6$.
 We denote by 
 \begin{eqnarray}
  \mathcal{C}_k^p({SM})=\left\{\varphi\in C^{\infty}({SM}): \forall i=0,1,...,k, \int_{SM}|\hat{\nabla}^i\varphi|^p\eta_F<\infty\right\}.
 \end{eqnarray}
 Because $SM$ is compact then $\mathcal{C}_k^p({SM})=C^{\infty}({SM})$
\begin{definition}\textbf{$($Sobolev spaces$)$}.
 Given the space $SM$ of all constant directions tangent on 3-Finsler manifold $M$, an integer $k\geq0$ and a real number $p\geq1$, 
 we define the Sobolev space $W^{k,p}(SM)$ as the completion of $\mathcal{C}_k^p({SM})$ with respect to the norm 
 $$||\varphi||_{W^{k,p}(SM)}=||\varphi||_{k,p}=\sum_{i=0}^k||\hat{\nabla}^i\varphi||_p.$$
 Set $C_c^{\infty}({SM})$ the space of $C^{\infty}$ functions with compact support in $SM$. 
 $W_0^{k,p}(SM)$ is the closure of $C_c^{\infty}({SM})$ in $W^{k,p}(SM)$.
\end{definition}
\subsection{Green's function of the Finslerian  Laplacian}
%By the remark (\ref{r}), we can (adapt?) the Green's function for a compact Finsler manifold $M$ when all directions tangent on $M$ are constant.
\begin{definition}\textbf{$($Green's function$)$}.
 Let $\overline{SM}^N$ be a $N:=2n-1$-dimensional 
 compact Riemannian manifold with boundary of classe $C^{\infty}$. The Green's function $G((x,[y]),(z,[y]))$ 
 of the horizontal Laplacian $\hat{\Delta}$ is the function which satisfies in $SM\times SM$:
 $$\hat{\Delta}_{(z,[y])}G((x,[y]),(z,[y]))=\delta_{(x,[y])}((z,[y]))$$ (where $\delta_{(x,[y])}$ 
 is the Dirac function at ${(x,[y])}$) and which vanishes on the boundary for $(x,[y])$ and $(z,[y])$ 
 belonging to $\partial(SM)$.
\end{definition}
\begin{proposition}\textbf{$($Aubin's theorem$)$}\cite{bi2019}.
Let $(M,F)$ be a $C^{\infty}$  compact Finsler manifold of dimension $n$ and $SM$ the space 
(of dimension $N:=2n-1$) of all constant directions tangent on $M$.
 %Let $SM^N$ be a $C^{\infty}$ Riemannian manifold. 
 There existes $G((x,[y]),(z,[y]))$ a Green's function of the Laplacian which has the following properties:
  \begin{itemize}
   \item [(1)] For all functions $\varphi\in C^2(SM)$, and for all $(x,[y])\in SM$:
   \begin{eqnarray}
    \varphi((x,[y]))&=&\frac{1}{vol(SM)}\int_{SM}\varphi((z,[y]))\eta_F((z,[y]))\nonumber\\
                    & &+\int_{SM}G((x,[y]),(z,[y]))\hat{\Delta}\varphi((z,[y]))\eta_F((z,[y])).
   \end{eqnarray}
\item [(2)] $G((x,[y]),(z,[y]))$ is $C^{\infty}$ on $SM\times SM \setminus \{(x,[y]),(z,[y]):(x,[y])=(z,[y])\}.$
\item [(3)] There exists a constant $c$ such that:\\
\noindent
$|G((x,[y]),(z,[y]))|=c(1+|log r|)$  for $N=2$  and\\
$|G((x,[y]),(z,[y]))|<c r^{2-N}$ for $N\geq 3,$ \\
$|\hat{\nabla}_{(z,[y])}G((x,[y]),(z,[y]))|<cr^{1-N},$\\
$|\hat{\nabla}_{(z,[y])}^2G((x,[y]),(z,[y]))|<cr^{-N},$ with $r=d((x,[y]),(z,[y])).$\\
\item [(4)] There exists a constant $A$ such that $G((x,[y]),(z,[y]))\geq A$. Because the Green's function is defined up 
to a constant, we can thus choose the Green's function every where positive.
\item [(5)] $\int_{SM}G((x,[y]),(z,[y]))\eta_F((z,[y]))=Constante$. 
We can choose the Green's function so that its integral equals zero.
\item [(6)] $G((x,[y]),(z,[y]))=G((z,[y]),(x,[y]))$. 
  \end{itemize}

\end{proposition}

       \begin{theorem}\label{tp}
        Let $(M^3,F)$ be a 3-dimensional closed Finsler manifold and $\hat{\textbf{S}}$ his constant nonpositive horizontal scalar curvature. 
        Let $\widetilde{a_i}, i=0,1,2$ be the heat invariants for the sphere bundle $SM$ of $M$.
        Suppose $\alpha_0,\alpha_1,\alpha_2$ and $\Lambda$ are positives constants, and $\lambda_1(\widetilde F)$ is 
        the first positive eigenvalue of the horizontal Laplacian $\hat{\Delta}$. If $\widetilde{F}=\varphi^2F$
        ($\varphi\in C^2(M), \varphi\neq0, \varphi\neq1$) is a 
        Finsler metric on $M$ satisfying
        \begin{itemize}
         \item[(i)] $\widetilde{a}_0=\alpha_0,$
         \item[(ii)] $\widetilde{a}_1 \leq \alpha_1,$
         \item[(iii)] $\widetilde{a}_2 \leq \alpha_2$
         \item[(iv)] $0<\Lambda \leq \lambda_1(\widetilde F),$
        \end{itemize}
        then there exist constants $c_1, c_2$ and $c_3$ 
        (depending only on $\alpha_0,\alpha_1,\alpha_2$ and $\Lambda$) such that 
        \begin{itemize}
         \item[(v)] $0 < c_1 \leq \varphi((x,[y]))\leq c_2~~~\forall~(x,[y])\in SM,$
         \item[(vi)]
%          \begin{eqnarray}
        $  ||\varphi||_{W^{2,2}(SM)} \leq c_3,$
%           \end{eqnarray} 
  where
           \begin{eqnarray}
           \label{1000}
           ||\varphi||_{W^{2,2}(SM)}
         =\left(\int_{SM}|\varphi|^2\eta_F+\int_{SM}|\hat{\nabla}\varphi|^2\eta_F
         +\int_{SM}|\hat{\nabla}^2\varphi|^2\eta_F
         \right)^{\frac{1}{2}}.
         \end{eqnarray} 
        \end{itemize}
       \end{theorem}
       Note that, in (\ref{1000}), the horizontal covariant derivative $\hat{\nabla}$ are understood to be with respect 
       to the Finsler metric $F$ and $|\varphi|=|\varphi((x,[y]))|~~\forall~(x,[y])\in SM$.
      %    \subsection{Some norm estimates}
     %Our goal in this paragraph is to bound the conformal factor $\varphi$ apriori in terms of the heat invariants 
     %$\widetilde{a_i}, i=0,1,2$ and the geometry of $SM$.
     \begin{lemma}
     Let $(M^3,F)$ be a 3-dimensional closed Finsler manifold and $SM$ its sphere bundle. \\
     Let $\widetilde{F}=\varphi^2F$ be a conformal Finsler metric of $F$.
     Then $\int_{SM}\varphi^2\eta_F$ and $\int_{SM}|\hat{\nabla}\varphi|^2\eta_F$ are apriori bounded.
     \end{lemma}
\begin{prof7}
 From the heat invariant $\widetilde{a_0}$, by applying H\"{o}lder inequality and by using that $1<\frac{2(N+1)}{N-3}$ we have 
 \begin{eqnarray}
  \int_{SM}\varphi^2\eta_F &\leq& \left(\int_{SM}|\varphi|^q\eta_F\right)^{\frac{2}{q}}
  \left(\int_{SM}1\eta_F\right)^{1-\frac{2}{q}},~~~ q=\frac{2(N+1)}{N-3}=6\nonumber\\
                 %    &=& \left(\int_{SM}|\varphi|^{\frac{2(N+1)}{N-3}}\eta_F\right)^{\frac{2}{{\frac{2(N+1)}{N-3}}}}
                  %   \left(\int_{SM}\eta_F\right)^{1-\frac{2}{{\frac{2(N+1)}{N-3}}}}\nonumber\\
                     &=& \widetilde{a_0}\times vol(SM)^{1-\frac{2}{{\frac{2(N+1)}{N-3}}}}\nonumber\\
                     &=& c
 \end{eqnarray}
where $c$ depends on $\widetilde{a_0}$ and the volume on $SM$.\\

We have just seen that $\hat{\widetilde{\textbf{S}}}$ satisfies to the equation (\ref{2}).
% $c(N)\hat{\Delta}\varphi+\hat{\textbf{S}}\varphi=\hat{\widetilde{\textbf{S}}}\varphi^{\frac{N+5}{N-3}}$. 
From the heat invariant\\ 
$\widetilde{a_1}=\frac{1}{6}\int_{SM}\hat{\widetilde{\textbf{S}}}\eta_{\widetilde F}$ we have
\begin{eqnarray}
 \widetilde{a_1}&=&\frac{1}{6}\int_{SM}(\hat{\textbf{S}}\varphi-c(N)\hat{\Delta}\varphi)\varphi^{-\frac{N+5}{N-3}}.\varphi^{
 \frac{2N+2}{N-3}}\eta_F,~~~~~~~c(N)=\frac{4(N-1)}{N-3}\nonumber\\
  %  &=&\frac{1}{6}\int_{SM}(\hat{\textbf{S}}\varphi-c(N)\hat{\Delta}\varphi)\varphi\eta_F\nonumber\\
    &=&\frac{1}{6}\int_{SM}(\hat{\textbf{S}}\varphi^2+c(N)|\hat{\nabla}\varphi|^2)\eta_F\nonumber\\
     &=&\frac{1}{6}\int_{SM}(\hat{\textbf{S}}\varphi^2+8|\hat{\nabla}\varphi|^2)\eta_F
\end{eqnarray}
when $n=3$. We obtain
\begin{eqnarray}
 \frac{4}{3}\int_{SM}|\hat{\nabla}\varphi|^2\eta_F&=&\widetilde{a_1}-\frac{1}{6}\int_{SM}\hat{\textbf{S}}\varphi^2\eta_F\nonumber\\
                                    &\leq& c
\end{eqnarray}
where $c$ is a constant depending on $\widetilde{a_1},\hat{\textbf{S}}$ and a bound of $\int_{SM}\varphi^2\eta_F$.\\

\noindent
 Hence $\int_{SM}\varphi^2\eta_F$ is bounded in terms of $\widetilde{a_0}, vol(SM)$
 while $\int_{SM}|\hat{\nabla}\varphi|^2\eta_F$ is bounded 
in terms of $\widetilde{a_1}, \hat{\textbf{S}}$ and a bound of $\int_{SM}\varphi^2\eta_F$.
\end{prof7}
\begin{lemma}
 Let $(M^3,F)$ be a 3-dimensional closed Finsler manifold and $SM$ its sphere bundle.
 Suppose $M^3$ has constant nonpositive horizontal scalar curvature. Then,\\
 $\int_{SM}(\hat{\widetilde{\textbf{S}}})^2\eta_{\widetilde F},~~~\int_{SM}\varphi^{-4}(\hat{\Delta}\varphi)^2\eta_F$~~~
 and~~~$\int_{SM}(\hat{\textbf{S}})^2\varphi^{-2}\eta_F$~~~are apriori bounded.
\end{lemma}
\begin{prof7}
 From the heat invariant $\widetilde{a_2}$, we have $\int_{SM}(\hat{\widetilde{\textbf{S}}})^2\eta_{\widetilde{F}}\leq 120\widetilde{a_2}$. But
 \begin{eqnarray}
  \int_{SM}(\hat{\widetilde{\textbf{S}}})^2\eta_{\widetilde{F}}&=&\int_{SM}|\hat{\textbf{S}}\varphi-c(N)\hat{\Delta}\varphi|^2\varphi^{\frac{-2(N+5)}{N-3}}
                                        \varphi^{\frac{2(N+1)}{N-3}}\eta_F\nonumber\\
                                    %&=&\int_{SM}|\hat{\textbf{S}}\varphi-8\hat{\Delta}\varphi|^2\varphi^{-10}
                                     %   \varphi^{6}\eta_F\nonumber\\
      %&=&\int_{SM}(\hat{\textbf{S}}\varphi)^2\varphi^{-4}\eta_F-16\int_{SM}\hat{\textbf{S}}\varphi^{-3}(\hat{\Delta}\varphi)\eta_F
      %+64\int_{SM}(\hat{\Delta}\varphi)^2\varphi^{-4}\eta_F\nonumber\\
       &=&\int_{SM}(\hat{\textbf{S}})^2\varphi^{-2}\eta_F-16\int_{SM}\hat{\textbf{S}}\varphi^{-3}(\hat{\Delta}\varphi)\eta_F
      +64\int_{SM}(\hat{\Delta}\varphi)^2\varphi^{-4}\eta_F.
 \end{eqnarray}
Using that $\hat{\textbf{S}}$ is constant and by integration by parts, we can rewrite the middle term as 
 \begin{eqnarray}
  I_2&=&-16\int_{SM}\hat{\textbf{S}}\varphi^{-3}(\hat{\Delta}\varphi)\eta_F\nonumber\\
     %&=&-16(-1)(-3)\hat{\textbf{S}}\int_{SM}|\hat{\nabla}\varphi|^2\varphi^{-4}\eta_F\nonumber\\
     &=&-48 \hat{\textbf{S}}\int_{SM}|\hat{\nabla}\varphi|^2\varphi^{-4}\eta_F.
 \end{eqnarray}
 When $\hat{\textbf{S}}=0$ then $I_2=0$ and when $\hat{\textbf{S}}<0$ then $I_2> 0$.
  We have 
  \begin{eqnarray}
  \int_{SM}(\hat{\textbf{S}})^2\varphi^{-2}\eta_F+64\int_{SM}(\hat{\Delta}\varphi)^2\varphi^{-4}\eta_F 
  \leq\int_{SM}(\hat{\widetilde{\textbf{S}}})^2\eta_{\widetilde{F}}\leq 120\widetilde{a_2}.
  \end{eqnarray}
  So $\int_{SM}(\hat{\widetilde{\textbf{S}}})^2\eta_{\widetilde F}$ and hence
%  \begin{eqnarray}
 $ \int_{SM}\varphi^{-4}(\hat{\Delta}\varphi)^2\eta_F$ and $\int_{SM}(\hat{\textbf{S}})^2\varphi^{-2}\eta_F$
%  \end{eqnarray}
   are apriori bounded, as desired.
\end{prof7}

Now we have
\begin{lemma}
 Let $(M^3,F)$ be a 3-closed Finsler manifold and $SM$ its sphere bundle. 
 When the horizontal scalar curvature $\hat{\textbf{S}}$ is negative $\varphi$ is pointwise apriori bounded away from 0.
\end{lemma}
\begin{prof7}
Let $\overline{(\frac{1}{\varphi})}=\frac{1}{vol(SM)}\int_{SM}\frac{1}{\varphi}\eta_F$ be the average value of 
 $\frac{1}{\varphi}$. Then when dimension $M=3$ and $\hat{\textbf{S}}<0$ we have 
  \begin{eqnarray}
 \overline{\left(\frac{1}{\varphi}\right)}\leq c\left(\int_{SM}\frac{1}{\varphi^2}\eta_F\right)^{\frac{1}{2}}
   \end{eqnarray}
   by Cauchy-Schwartz, so that 
 $\overline{(\frac{1}{\varphi})}$ is apriori bounded, by the last peceding Lemma.\\
 
 We now compute, for any point $(x,[y])\in SM$, 
\begin{eqnarray}
\frac{1}{\varphi((x,[y]))}=\overline{\left(\frac{1}{\varphi}\right)}+\int_{SM}G\left((x,[y]),(z,[y])\right)
 \left[-\hat{\Delta}\left(\frac{1}{\varphi((z,[y]))}\right)\right]\eta_F((z,[y])) 
\end{eqnarray}
 where $G\left((x,[y]),(z,[y])\right)$ is the Green's function of the horizontal Laplacian $\hat{\Delta}$ on $SM$. Since
 \begin{eqnarray}
  \hat{\Delta}\left(\frac{1}{\varphi}\right)&=&-\frac{1}{\varphi^2}\hat{\Delta}\varphi+\frac{2|\hat{\nabla}\varphi|^2}{\varphi^3}\nonumber\\
                       &=&-\frac{1}{8\varphi^2}\left(\hat{\textbf{S}}\varphi-\hat{\widetilde{\textbf{S}}}\varphi^5\right)
                       +\frac{2|\hat{\nabla}\varphi|^2}{\varphi^3}
 \end{eqnarray}
by Yamabe equation. We get from Green's fonction that for any point $(x,[y])\in SM$,
\begin{eqnarray}
 \frac{1}{\varphi((x,[y]))}-\overline{\left(\frac{1}{\varphi}\right)}
 &=&\int_{SM}G((x,[y]),(z,[y])) \nonumber\\
 & &\times\left[\frac{1}{8}\hat{\textbf{S}}(\varphi)^{-1}
  -\frac{1}{8}\hat{\widetilde{\textbf{S}}}\varphi^3
  -\frac{2|\hat{\nabla}\varphi|^2}{\varphi^3}\right]\eta_F((z,[y]))\nonumber\\
  %&\leq&\frac{1}{8}\int_{SM}G((x,[y]),(z,[y]))\hat{\textbf{S}}\varphi^{-1}\eta_F((z,[y]))\nonumber\\
  %& &-\frac{1}{8}\int_{SM}G((x,[y]),(z,[y]))\hat{\widetilde{\textbf{S}}}\varphi^3\eta_F((z,[y]))\nonumber\\
  &\leq &\frac{1}{8}\left(\int_{SM}(G((x,[y]),(z,[y])))^2\eta_F((z,[y]))\right)^\frac{1}{2}\nonumber\\
 & & \times\left[\left(\int_{SM}(\hat{\textbf{S}})^2\varphi^{-2}\eta_F((z,[y]))\right)^\frac{1}{2}\right.\nonumber\\
 & & +\left.\left(\int_{SM}(\hat{\widetilde{\textbf{S}}})^2\varphi^6\eta_F((z,[y]))\right)^{\frac{1}{2}}\right],
\end{eqnarray}
where the two terms in the bracket are bounded by Lemma $5.2$, and nothing that, for dimension of $M=3$, $G((x,[y]),(z,[y]))$ is
square integrable. 
It follows that there exists a constant $c_1>0$ which does not depend on 
such that for any $(x,[y])\in SM$ $$0<c_1\leq \varphi((x,[y])).$$
\end{prof7}

\subsection{Proof of Theorem 5.1}
We proceed in two steps in order to find an upper bound for $\varphi$
\begin{lemma}
 Let $(M^3,F)$ be a 3-dimensional closed Finsler manifold and $SM$ its sphere bundle. Assume $\varphi$ be a positive function 
 on $SM$ which satisfies the hypothesis of theorem $5.1$, 
 there exists some $\varepsilon_0 > 0$ sufficiently small and a constant $c(a_0, a_1, a_2, \lambda_1)$ so that 
\begin{eqnarray}
 \int_{SM}\varphi^{6+\varepsilon_0}\eta_F\leq c.
 \end{eqnarray}
\end{lemma}
\begin{prof7}
 Let $\phi=\varphi^{1+\varepsilon}$ with $\varepsilon$ to be determine later. We have, from Sobolev inequality for $\phi$ \cite{bi11}, that 
 \begin{eqnarray}
 \label{98}
\left(\int_{SM}\phi^6\eta_F\right)^{\frac{1}{3}}\leq C_1\int_{SM}|\hat{\nabla}\phi|^2\eta_F+C_2\int_{SM}\phi^2\eta_F,
 \end{eqnarray}
where $C_1, C_2$ depend only on the geometry of $SM$.
On the other hand, multiplying the Yamabe (type) equation $8\hat{\Delta}\varphi-\hat{\textbf{S}}\varphi=-\hat{\widetilde{\textbf{S}}}\varphi^5$ by 
$\varphi^{1+2\varepsilon}$, we have 
\begin{eqnarray}
\label{99}
 8\varphi^{1+2\varepsilon}.\hat{\Delta}\varphi-\hat{\textbf{S}}\varphi^{2+2\varepsilon}=-\hat{\widetilde{\textbf{S}}}\varphi^4.\varphi^{2+2\varepsilon}.
\end{eqnarray}
By integration of (\ref{99}) and by using integration by parts, we get
\begin{eqnarray}
\label{100}
 8\int_{SM}\varphi^{1+2\varepsilon}.\hat{\Delta}\varphi\eta_F &=& 8\int_{SM}\hat{\Delta}\left(\phi^{\frac{1}{1+\varepsilon}}\right)
       \phi^{\frac{1+2\varepsilon}{1+\varepsilon}}\eta_F\nonumber\\
   % &=& -8\int_{SM}\left(\frac{1}{1+\varepsilon}\phi^{\frac{-\varepsilon}{1+\varepsilon}}\hat{\nabla}\phi\right)
         %\left(\frac{1+2\varepsilon}{1+\varepsilon}\right)\phi^{\frac{\varepsilon}{1+\varepsilon}}\hat{\nabla}\phi\eta_F\nonumber\\
   &=& -\frac{8(1+2\varepsilon)}{(1+\varepsilon)^2}\int_{SM}\left|\hat{\nabla}\phi\right|^2\eta_F.\
\end{eqnarray}
% Combining (\ref{100}) and (\ref{99}), we find 
so that
\begin{eqnarray}
 \label{101}
 & & \frac{8(1+2\varepsilon)}{(1+\varepsilon)^2}\int_{SM}\left|\hat{\nabla}\phi\right|^2\eta_F+\hat{\textbf{S}}\int_{SM}\phi^2\eta_F
 =\int_{SM}\hat{\widetilde{\textbf{S}}}\varphi^4\phi^2\eta_F\nonumber\\
 \end{eqnarray}
 or equivalently,
\begin{eqnarray}
  \int_{SM}\left|\hat{\nabla}\phi\right|^2\eta_F = \frac{(1+\varepsilon)^2}{8(1+2\varepsilon)}
 \int_{SM}\hat{\widetilde{\textbf{S}}}\varphi^4\phi^2\eta_F-\frac{(1+\varepsilon)^2}{8(1+2\varepsilon)}\hat{\textbf{S}}\int_{SM}\phi^2\eta_F.\nonumber\\
\end{eqnarray}
(\ref{101}) in (\ref{98}) 
\begin{eqnarray}
\label{102}
 \left(\int_{SM}\phi^6\eta_F\right)^{\frac{1}{3}}&\leq& C_1\left[\frac{(1+\varepsilon)^2}{8(1+2\varepsilon)}
 \int_{SM}\hat{\widetilde{\textbf{S}}}\varphi^4\phi^2\eta_F-\frac{(1+\varepsilon)^2}{8(1+2\varepsilon)}\hat{\textbf{S}}\int_{SM}\phi^2\eta_F\right]\nonumber\\
 & &+C_2\int_{SM}\phi^2\eta_F.
\end{eqnarray}
For $\varepsilon < 2$, $\int_{SM}\phi^2\eta_F$ is bounded by %$C(\varepsilon)\int_{SM}\varphi^6\eta_F,$ (Lemma $5.4)$, 
some constant depending on $\varepsilon, \widetilde{a_0}$ and the volume of $SM$. Hence 
\begin{eqnarray}
 \left(\int_{SM}\phi^6\eta_F\right)^{\frac{1}{3}}&\leq& 
 C_1\frac{(1+\varepsilon)^2}{8(1+2\varepsilon)}\int_{SM}\hat{\widetilde{\textbf{S}}}\varphi^4\phi^2\eta_F+C_2(\varepsilon).\
\end{eqnarray}

We now seek to bound $\int_{SM}\hat{\widetilde{\textbf{S}}}\varphi^4\phi^2\eta_F$ from above in terms of $\int_{SM}\phi^6\eta_F$, and to 
that end we consider a large constant $b$ to be chosen later, and we consider the sets where $|\hat{\widetilde{\textbf{S}}}|\geq 0$ and 
$|\hat{\widetilde{\textbf{S}}}|<0$.\\
For $|\hat{\widetilde{\textbf{S}}}|\geq 0$,
\begin{eqnarray}
 \int_{|\hat{\widetilde{\textbf{S}}}|\geq b}\hat{\widetilde{\textbf{S}}}\varphi^4\phi^2\eta_F &\leq& \left(
 \int_{|\hat{\widetilde{\textbf{S}}}|\geq b}(\hat{\widetilde{\textbf{S}}})^2\varphi^6\eta_F\right)^{\frac{1}{2}}\left(
 \int_{|\hat{\widetilde{\textbf{S}}}|\geq b}\varphi^6\eta_F\right)^{\frac{1}{6}}\left(
 \int_{|\hat{\widetilde{\textbf{S}}}|\geq b}\phi^6\eta_F\right)^{\frac{1}{3}}\nonumber\\
 &\leq &\left(a_2\right)^{\frac{1}{2}}\left(\frac{a_2}{b^22}\right)^{\frac{1}{6}}
 \left(\int_{SM}\phi^6\eta_F\right)^{\frac{1}{3}},\
\end{eqnarray}
while for $|\hat{\widetilde{\textbf{S}}}|<b$ we have 
\begin{eqnarray}
 \int_{|\hat{\widetilde{\textbf{S}}}|<b}\hat{\widetilde{\textbf{S}}}\varphi^4\phi^2\eta_F &\leq& b\int_{SM}\varphi^4\phi^2\eta_F,
\end{eqnarray}
so that 
\begin{eqnarray}
\label{165}
 \int_{SM}\hat{\widetilde{\textbf{S}}}\varphi^4\phi^2\eta_F &\leq&\left(a_2\right)^{\frac{1}{2}}\left(\frac{a_2}{b_2}\right)^{\frac{1}{6}}
 \left(\int_{SM}\phi^6\eta_F\right)^{\frac{1}{3}}+b\int_{SM}\varphi^4\phi^2\eta_F.\
\end{eqnarray}
Let $\lambda_1=\lambda_1(\hat{\widetilde{\Delta}})\equiv \Lambda$ denote 
the first eigenvalue of $\hat{\Delta}$ acting on functions, then to estimate
 $\int_{SM}\varphi^4\phi^2\eta_F$, we apply the Rayleigh-Ritz characterization of $\lambda_1$ \cite{}:
 \begin{eqnarray}
  \lambda_1(\hat{\widetilde{\Delta}})\leq \frac{\int_{SM}|\hat{\nabla}_{\widetilde F}\psi|^2\eta_{\widetilde F}}
  {\int_{SM}\psi^2\eta_{\widetilde F}-(vol(SM))^{-1}\left(\int_{SM}\psi\eta_{\widetilde F}\right)^2}
 \end{eqnarray}
or equivalently, 
\begin{eqnarray}
 \int_{SM}\psi^2\eta_{\widetilde F}\leq (vol(SM))^{-1}\left(\int_{SM}\psi\eta_{\widetilde F}\right)^2
 +\frac{1}{\lambda_1}\int_{SM}|\hat{\nabla}_{\widetilde F}\psi|^2\eta_{\widetilde F}
\end{eqnarray}
to $\psi=\varphi^{\varepsilon}$ to obtain
\begin{eqnarray}
 \int_{SM}\varphi^4\phi^2\eta_F\leq \left(\int_{SM}\varphi^6\eta_F\right)^{-1}
 \left(\int_{SM}\varphi^{6+\varepsilon}\eta_F\right)^{2}+
 \frac{1}{\lambda_1}\int_{SM}|\hat{\nabla}\varphi^{\varepsilon}|^2\eta_{\widetilde F}.
\end{eqnarray}
But 
\begin{eqnarray}
 \int_{SM}|\hat{\nabla}_{\widetilde F}\varphi^\varepsilon|^2\eta_{\widetilde F}&=&
 \int_{SM}|\hat{\nabla}\varphi^\varepsilon|^2\varphi^2\eta_F\nonumber\\
% &=&\frac{\varepsilon^2}{(1+\varepsilon)^2}\int_{SM}|\hat{\nabla}\varphi^{1+\varepsilon}|^2\eta_F\nonumber\\
 &=&\frac{\varepsilon^2}{8(1+\varepsilon)}\left(\int_{SM}\hat{\widetilde{\textbf{S}}}\varphi^4\phi^2\eta_F-
 \hat{\textbf{S}}\int_{SM}\phi^2\eta_F\right)\nonumber\\
\end{eqnarray}
where again we have used $8\hat{\Delta}\varphi-\hat{\textbf{S}}\varphi=-\hat{\widetilde{\textbf{S}}}\varphi^5$ to obtain the last line.\\
To estimate $\int_{SM}\varphi^{6+\varepsilon}\eta_F$, we choose $c$ as in Lemma $5.4$, with $\varphi-c\geq 0$, so that
\begin{eqnarray}
 \int_{SM}\varphi^{6+\varepsilon}\eta_F &=& \int_{SM}(\varphi^6-c^6)\varphi^\varepsilon\eta_F
                                  +\int_{SM}c^6\varphi^\varepsilon\eta_F \nonumber\\
 &\leq& \left[\int_{SM}(\varphi^6-c^6)\varphi^{2\varepsilon}\eta_F\right]^{\frac{1}{2}}
 \left[(\varphi^6-c^6)\eta_F\right]^{\frac{1}{2}}+\int_{SM}c^6\varphi^\varepsilon\eta_F,
\end{eqnarray}
where the inequality comes from Cauchy-Schwartz and the positivity of $(\varphi^6-c^6)$. 
Squaring and using the inequality $2AB\leq\delta A^2+(\frac{1}{\delta})B^2$, we obtain
\begin{eqnarray}
\left(\int_{SM}\varphi^{6+\varepsilon}\right)^{\frac{1}{2}}\leq (1+\delta)\left[
\int_{SM}(\varphi^6-c^6)\varphi^{2\varepsilon}\eta_F\right]\left[(\varphi^6-c^6)\eta_F\right]
+(1+\delta)\left[\int_{SM}c^6\varphi^\varepsilon\eta_F\right]
\end{eqnarray}
where $\delta$ will be chosen late. But $\int_{SM}(\varphi^6-c^6)\eta_F=\alpha.\int_{SM}\varphi^6\eta_F$ 
where $\alpha =1-\frac{c^6.vol_{F}(SM)}{vol_{\widetilde{F}}(SM)}$ is a positive constant less than 1 and we conclude 
\begin{eqnarray}
 \left(\int_{SM}\varphi^6\eta_F\right)^{-1}\left(\int_{SM}\varphi^{6+\varepsilon}\eta_F\right)^2
 &\leq& (1+\delta).\alpha\left(\int_{SM}\varphi^{6+2\varepsilon}\eta_F\right)+(1+\frac{1}{\delta}).(constant).
\end{eqnarray}
Choosing $\delta$ so that $(1+\delta)\alpha=(1-\nu) < 1$, we then have 
\begin{eqnarray}
\int_{SM}\varphi^4\phi^2\eta_F\leq(1-\nu)\int_{SM}\varphi^4\phi^2\eta_F+\frac{\varepsilon^2}{8(1+\varepsilon)}
\int_{SM}\hat{\widetilde{\textbf{S}}}\varphi^4\phi^2\eta_F+(constant),
\end{eqnarray}
so that
\begin{eqnarray}
\label{200}
 \int_{SM}\varphi^4\phi^2\eta_F\leq\frac{1}{\nu}\left[\frac{\varepsilon^2}{8(1+\varepsilon)}
\int_{SM}\hat{\widetilde{\textbf{S}}}\varphi^4\phi^2\eta_F+(constant)\right].
\end{eqnarray}
Combining (\ref{165}) and (\ref{200}), we see that 
\begin{eqnarray}
 \int_{SM}\hat{\widetilde{\textbf{S}}}\varphi^4\phi^2\eta_F\leq (a_2)^{\frac{1}{2}}\left(\frac{a_2}{b_2}\right)^{\frac{1}{6}}
 \left(\int_{SM}\phi^6\eta_F\right)^{\frac{1}{3}}+\frac{b}{\nu}\left[\frac{\varepsilon^2}{8(1+\varepsilon)}
\int_{SM}\hat{\widetilde{\textbf{S}}}\varphi^4\phi^2\eta_F+(constant)\right].
\end{eqnarray}
Choosing $b$ large so that the coefficient of $\left(\int_{SM}\phi^6\eta_F\right)^{\frac{1}{3}}$ is $<\varphi_1$ 
and then choosing $\varepsilon$ small so that the coefficient of $\int_{SM}\hat{\widetilde{\textbf{S}}}\varphi^4\phi^2\eta_F$ 
is $<\varphi_2$, we have
\begin{eqnarray}
 (1-\varphi_2)\int_{SM}\hat{\widetilde{\textbf{S}}}\varphi^4\phi^2\eta_F\leq
 \varphi_1\left(\int_{SM}\phi^6\eta_F\right)^{\frac{1}{3}}+(constant),
\end{eqnarray}
and comparing with the lower bound (\ref{102}), we have an upper bound for $\left(\int_{SM}\phi^6\eta_F\right)^{\frac{1}{3}}$
 as soon as $\frac{8(1+2\varepsilon)}{C_1(1+\varepsilon)^2}\geq\frac{\varphi_1}{1-\varphi_2}$,
 which achieved by taking $\varepsilon$ sufficiently small. This proves Lemma 4.
\end{prof7}

\begin{lemma}
 Let $(M^3,F)$ be a 3-dimensional closed Finsler manifold and $SM$ its sphere bundle. 
 There exists a positive constant $c_2(a_0, a_1, a_2, \lambda_1, \int_{SM}\varphi^{6+\varepsilon_0}\eta_F)$, such that
 $$\varphi\leq c_2$$ on $SM$.
\end{lemma}
\begin{prof7}
 Apply Green's function to equation (\ref{2}), we have, for all $(z,[y]) \in SM$
 \begin{eqnarray}
 \label{50}
  \varphi((x,[y]))&=&\overline{\varphi}-\int_{SM}\hat{\Delta}\varphi((z,[y]))
  G((x,[y]),(z,[y]))\eta_F((z,[y]))\nonumber\\
  &=&\overline{\varphi}+\frac{1}{8}\int_{SM}(\hat{\widetilde{\textbf{S}}}\varphi^5-\hat{\textbf{S}})G((x,[y]),(z,[y]))\eta_F((z,[y])).
 \end{eqnarray}
Since $\overline{\varphi}$ and $\int_{SM}\varphi G((x,[y]),(z,[y]))\eta_F((z,[y]))$ are apriori bounded, 
to bound $\varphi((x,[y]))$ it suffices to bound 
$\int_{SM}(\hat{\widetilde{\textbf{S}}}\varphi^5)G((x,[y]),(z,[y]))\eta_F((z,[y]))$.
By H\"{o}lder, we have 
\begin{eqnarray}
\label{51}
\int_{SM}(\hat{\widetilde{\textbf{S}}}\varphi^5)G((x,[y]),(z,[y]))\eta_F((z,[y]))
 &\leq&\left(\int_{SM}|\hat{\widetilde{\textbf{S}}}\varphi^5|^{p}\eta_F\right)^{\frac{1}{p}}\times \nonumber\\
 & &\left(\int_{SM}|G((x,[y]),(z,[y]))|^{q}\eta_F((z,[y]))\right)^{\frac{1}{q}}\nonumber\\
 &=&||\hat{\widetilde{\textbf{S}}}\varphi^5||_{p}||G((x,[y]),.)||_{q}
\end{eqnarray}
for $\frac{1}{p}+\frac{1}{q}=1.$ Choose $p=\frac{3}{5}+\delta, \delta=\frac{\varepsilon}{16+2\varepsilon}$ 
we found $q<\frac{5}{2}$ so that $||G((x,[y]),.)||_{q}$ is finite for $q$, so we must bound 
$||G((x,[y]),.)||_{q}$ for $q>\frac{5}{2}$. \\
It well known that
$|G((x,[y]),(z,[y_2]))|\leq C / d^3((x,[y]),(z,[y_2]))$ for some constant $C$. Recall
the following estimate: for $h(z)=\int_{\mathbb{R}^5}(f(x)/||x-z||)dz$ we have 
\begin{eqnarray}
 ||h||_r\leq c(p)||f||_{p},~~~r=d((x,[y]),(z,[y_2]))
\end{eqnarray}
when $\frac{1}{r}=\frac{1}{p}-\frac{3}{5}$ with $r>1$.
Since $G((x,[y]),(z,[y_2]))\leq C_2(r)$, we have 
\begin{eqnarray}
 ||\varphi||_r&\leq& C(r)+||\hat{\widetilde{\textbf{S}}}\varphi^5||_{p}||G((x,[y]),(z,[y]))||_{q}\nonumber\\
               &\leq& C(r)+C((z,[y]))||\hat{\widetilde{\textbf{S}}}\varphi^5||_p\nonumber\\
               %&\leq& C(r)+C((z,[y]))||\hat{\widetilde{\textbf{S}}}\varphi^5||_p^p\nonumber\\
               &=&C(r)+C((z,[y])).\int_{SM}(\hat{\widetilde{\textbf{S}}}\varphi^5)^p\eta_F\nonumber\\
               %&\leq&C(r)+C((z,[y])).\left(\int_{SM}(\hat{\widetilde{\textbf{S}}})^2\varphi^6\eta_F\right)^{p/2}.
               %\left(\int_{SM}\varphi^r\eta_F\right)^{1-p/2}\nonumber\\
               &=&C(r)+C((z,[y])).(A_2)^{p/2}.\left(\int_{SM}\varphi^r\eta_F\right)^{1-p/2}\
\end{eqnarray}
for $r=\frac{4p}{2-p}$ with $r>1$. We see that a bound for $||\varphi||_r$ implies a bound for $||\hat{\widetilde{\textbf{S}}}\varphi^5||_{p}$
for $$p=\frac{2r}{4+r}$$ and the a bound for $$\widetilde{r}=\frac{6r}{12-r}.$$ 
Starting this procedure with $r_0=6+\varepsilon$, with the bound on $||\varphi||_{r_0}$, given by Lemma $5.4$, we set 
$$r_{k+1}=\frac{6r_k}{12-r_k}$$ and $$p_{k}=\frac{2r_k}{4+r_k}.$$
Nothing that 
\begin{eqnarray}
 r_{k+1}-r_k=\frac{6r_k}{12-r_k}-r_k=r_k\frac{r_k-6}{12-r_k},
\end{eqnarray}
as long as $r_k<12$, the $r_k$ increase at least geometrically. 
It follows that there is a $k_0$, depending on $\varepsilon$, so that $r_0>12$, and hence $q_{k_0}>3/2$. 
For this $k_0$, (\ref{50}) and (\ref{51}) provide the $C_2$ required by Lemma. This finish the proof of Lemma $5.5$
\end{prof7}

End of the proof of theorem $5.3$. From the Lemma $5.3$ and $5.4$ we have 
$$0<c_1\leq\varphi((x,[y]))\leq c_2~~~\forall~(x,[y])\in SM.$$ 
From 
\begin{eqnarray}
 \alpha_2\geq \widetilde{a_2}\geq(\hat{\widetilde{\textbf{S}}})^2\varphi^6\eta_F
 =\int_{SM}\left(\frac{64(\hat{\Delta} \varphi)^2}{\varphi^4}-\frac{16\hat{\textbf{S}}(\hat{\Delta} \varphi)^3}{\varphi^4}+\frac{(\hat{\textbf{S}})^2}{\varphi^2}\right)\eta_F,
\end{eqnarray}
we conclude that 
\begin{eqnarray}
 \int_{SM}\frac{64(\Delta \varphi)^2}{\varphi^4}\eta_F\leq constant.
\end{eqnarray}
This, together with the uniform upper bound for $\varphi$, yields a bound for 
$\int_{SM}\Delta \varphi)^2\eta_F$. Thus $\int_{SM}[(\Delta \varphi)^2+\varphi^2]\eta_F$ is bounded; 
hence we have bound for $||\varphi||_{W^{2,2}(SM)}.$

\begin{proposition}
Let $(M^3,F)$ be a 3-dimensional closed Finsler manifold and $SM$ its sphere bundle of constant directions tangent on $M$.
%       For a compact Finsler manifold $(M^3,F)$ without boundary, 
Then the set of conformal metrics $\widetilde{F}_i=\varphi_i^2F$ with $\varphi_i$ satisfying 
 $(v), (vi)$ of theorem (\ref{tp}) is compact and the conditions $a_k(F)\leq \alpha_k, k=0,1,2$ forms a compact set in $C^2$ topology.
\end{proposition}

%\begin{remark}
 %This result was proved directly using the fist three heat invariants ($a_0,a_1$ and $a_2$) of the Finsler metric $F$
 %when all directions tangent $SM$ on $3$-closed manifold are constants. In our next work, we will prove the general case,
  %that is for every $a_k$.
%\end{remark}

\section*{Conflicts of Interest}
The authors declare no conflicts of interest regarding the publication of this work.

 \end{document}